\documentclass{amsart}
\usepackage{amssymb}
\usepackage{amssymb,amscd,amsthm}
\usepackage{latexsym}
\date{\today}

\newcommand{\Z}{{\mathbb Z}}
\newcommand{\R}{{\mathbb R}}
\newcommand{\C}{{\mathbb C}}
\newcommand{\N}{{\mathbb N}}

\newtheorem{theorem}{Theorem} [section]

\newtheorem{lemma}{Lemma}   [section]
\newtheorem{prop}{Proposition} [section]

\newtheorem{definition}{Definition} [section]
\begin{document}
\title[Limit-Periodic Schr\"odinger Operators]{Spectral Properties of Limit-Periodic\\Schr\"odinger Operators}

\author[D.\ Damanik]{David Damanik}

\address{Department of Mathematics, Rice University, Houston, TX~77005, USA}

\email{damanik@rice.edu}

\urladdr{www.ruf.rice.edu/$\sim$dtd3}

\author[Z.\ Gan]{Zheng Gan}

\address{Department of Mathematics, Rice University, Houston, TX~77005, USA}

\email{zheng.gan@rice.edu}

\urladdr{math.rice.edu/$\sim$zg2}

\thanks{D.\ D.\ and Z.\ G.\ were supported in part by NSF grant
DMS--0800100.}

\begin{abstract}
We investigate the spectral properties of Schr\"odinger operators
in $\ell2(\Z)$ with limit-periodic potentials. The perspective we
take was recently proposed by Avila and is based on regarding such
potentials as generated by continuous sampling along the orbits of
a minimal translation of a Cantor group. This point of view allows
one to separate the base dynamics and the sampling function. We
show that for any such base dynamics, the spectrum is of positive
Lebesgue measure and purely absolutely continuous for a dense set
of sampling functions, and it is of zero Lebesgue measure and
purely singular continuous for a dense $G_\delta$ set of sampling
functions.
\end{abstract}

\maketitle

\section{Introduction}

Schr\"odinger operators $H_\omega$ acting in $\ell2(\Z)$ with
dynamically defined potentials $V_\omega$ are given by
\begin{equation}\label{oper}
[H_\omega \psi](n) = \psi(n+1) + \psi(n-1) + V_\omega(n) \psi(n),
\end{equation}
where
\begin{equation}\label{pot}
V_\omega (n) = f(T^n \omega) , \quad \omega \in \Omega, \; n \in
\Z
\end{equation}
with a homeomorphism $T$ of a compact space $\Omega$ and a
continuous sampling function $f : \Omega \to \R$. It is often
beneficial to study the operators $\{H_\omega\}_{\omega \in
\Omega}$ as a family, as opposed to a collection of individual
operators. This is especially true if a $T$-ergodic probability
measure $\mu$ is chosen since the spectrum and the spectral type
of $H_\omega$ are always $\mu$-almost surely independent of
$\omega$ due to ergodicity. Alternatively, in the topological
setting, if $T$ is minimal and $f$ is continuous, then the
spectrum and the absolutely continuous spectrum are independent of
$\omega$; however, the point spectrum and the singular continuous
spectrum are in general not independent of $\omega$.

In this paper we will study the class of limit-periodic
potentials. As pointed out recently by Avila \cite{a}, these
potentials are defined by continuous sampling along the orbits of
a minimal translation of a Cantor group. We give the definitions
of these terms now and defer the discussion of their relation to
limit-periodic potentials to Section~\ref{s.profinite}.

Indeed, promoting Avila's approach to limit-periodic potentials is
part of our motivation for writing this paper because we believe
that it has the potential to produce further insight into
limit-periodic Schr\"odinger operators, in addition to the
beautiful results obtained in \cite{a}. Some developments in this
direction can be found in \cite{dg}.

\begin{definition}
We say that $\Omega$ is a Cantor group if it is a totally
disconnected compact Abelian topological group with no isolated
points. A map $T : \Omega \to \Omega$ is called a translation if
it is of the form $T \omega = \omega + \alpha$ for some $\alpha
\in \Omega$. A translation is called minimal if the orbit $\{ T^n
\omega : n \in \Z \}$ of every $\omega \in \Omega$ is dense in
$\Omega$.
\end{definition}

Given a minimal translation $T$ of a Cantor group $\Omega$ and a
continuous sampling function $f$, we consider the potentials
generated by \eqref{pot} and the associated operators given by
\eqref{oper}. A major advantage of Avila's approach to
limit-periodic potentials is that one can keep the base dynamics
$T : \Omega \to \Omega$ fixed and vary the sampling function. This
enables one to exhibit spectral properties that occur frequently,
that is, for a dense set of sampling functions, or even a dense
$G_\delta$ set of sampling functions. The following theorem
describes such spectral properties that occur frequently; it makes
statements about the Lebesgue measure of the spectrum, the
structure of the gaps of the spectrum, and the type of the
spectral measures.

\begin{theorem}\label{t.main}
Suppose $\Omega$ is a Cantor group and $T : \Omega \to \Omega$ is
a minimal translation. \\[1mm]
{\rm (a)} For a dense set of $f \in C(\Omega,\R)$ and every
$\omega \in \Omega$, the spectrum of $H_\omega$ is a Cantor set of
positive Lebesgue measure and $H_\omega$ has purely absolutely
continuous spectrum. \\
{\rm (b)} For a dense $G_\delta$-set of $f \in C(\Omega,\R)$ and
every $\omega \in \Omega$, the spectrum of $H_\omega$ is a Cantor
set of zero Lebesgue measure and $H_\omega$ has purely singular
continuous spectrum.
\end{theorem}

\noindent\textit{Remarks.} (i) A Cantor set is by definition a
closed set with empty interior and no isolated points. By general
principles, the only property that needs to be addressed in the
proof is the empty interior. 
\\[1mm]
(ii) The statement in (b) can be strengthened as follows. If
instead of $f$, one considers the one-parameter family $\{ \lambda
f \}_{\lambda > 0}$, then the same conclusion holds for the family
uniformly in $\lambda$. This is inferred easily from the proof.
\\[1mm]
(iii) The theorem above lists both known and new results for the
sake of completeness. Let us explain which aspects are new. Part
(a) is a version of results established in somewhat different form
by various authors for continuum Schr\"odinger operators in the
1980's; compare \cite{as,c1,m}. It is intermediate between the existing denseness
statements (again in the continuum setting) for all limit-periodic
potentials with variation in $\Omega$ allowed and for specific
families of limit-periodic potentials, where $\Omega$ is often
fixed but the class of considered $f$'s is restricted. Our proof
is close in spirit to the work of Avron and Simon \cite{as}. Most
claims in part (b) are due to Avila \cite{a}. What we add here is
the generic absence of eigenvalues.
\\[1mm]
(iv) It would clearly be of interest to determine how often the
third basic spectral type, pure point spectrum, occurs for
limit-periodic potentials. Some examples were constructed by P\"oschel \cite{p}.
Moreover, further results in this direction are stated in \cite{mc}, and again in \cite{c}.
We plan to explore limit-periodic Schr\"odinger operators with pure point spectrum
(both with and without exponentially decaying eigenfunctions) in a
future work.
\\[1mm]
(v) In the families $\{H_\omega\}_{\omega \in \Omega}$ covered by
Theorem~\ref{t.main}, the spectral type of $H_\omega$ is
independent of $\omega$. It is known that such uniformity may fail
in the more general class of almost periodic Schr\"odinger
operators. For example, if one takes $\Omega = \R / \Z$, $T \omega
\mapsto \omega + \alpha$ with a Diophantine number $\alpha$, and
$f(\omega) = 3 \cos (2 \pi \omega)$, then the associated operator
$H_\omega$ has pure point spectrum for Lebesgue almost every
$\omega$ and purely singular continuous spectrum for a dense
$G_\delta$ set of $\omega$'s. Thus, another interesting open
problem is the following: can a limit-periodic family
$\{H_\omega\}_{\omega \in \Omega}$ ever exhibit such a phenomenon
or, on the contrary, is the spectral type of limit-periodic
Schr\"odinger operators always independent of $\omega$?
\\[1mm]
(vi) All families $\{H_\omega\}_{\omega \in \Omega}$ covered by Theorem~\ref{t.main} have Cantor spectrum and one may again ask whether this is always the case or not. In fact, it is known that this is not always the case. Again in P\"oschel's paper \cite{p}, one may find examples of limit-periodic Schr\"odinger operators whose spectra have no gaps.
\\[1mm]
(vii) The present paper and \cite{dg} treat limit-periodic Schr\"odinger
operators in two complementary regimes. The operators discussed
here have zero Lyapunov exponents, while the operators discussed
in \cite{dg} have positive Lyapunov exponents. The main result of \cite{dg} shows that for every Cantor group $\Omega$ and every minimal translation $T$ on $\Omega$, there exists a dense set $\mathcal{F} \subset
C(\Omega,\R)$ such that for every $f \in \mathcal{F}$ and every $\lambda \in \R \setminus \{0\}$, the following statements hold true for the operators $H_\omega$ with potentials $V_\omega (n) = \lambda f(T^n \omega)$: the spectrum of $H_\omega$ has zero Hausdorff dimension and all spectral measures are purely singular continuous (for every $\omega \in \Omega$), and the Lyapunov exponent is a positive continuous function of the energy.

\section{Minimal Translations of Cantor Groups and the Description of Limit-Periodic
Potentials}\label{s.profinite}

In this section we recall how the one-to-one correspondence
between hulls of limit-periodic sequences and potential families
generated by minimal translations of Cantor groups and continuous
sampling functions exhibited by Avila in \cite{a} arises.

\begin{definition}
Let $S : \ell^\infty(\Z) \to \ell^\infty(\Z)$ be the shift
operator, $(S V)(n) = V(n+1)$. A two-sided sequence $V \in
\ell^\infty(\Z)$ is called periodic if its $S$-orbit is finite and
it is called limit-periodic if it belongs to the closure of the
set of periodic sequences. If $V$ is limit-periodic, the closure
of its $S$-orbit is called the hull and denoted by
$\mathrm{hull}_V$.
\end{definition}

The first lemma (see \cite[Lemma~2.1]{a}) shows how one can write
the elements of the hull of a limit-periodic function in the form
\eqref{pot} with a minimal translation $T$ of a Cantor group and a
sampling function $f \in C(\Omega,\R)$:

\begin{lemma}
Suppose $V$ is limit-periodic. Then, $\Omega : = \mathrm{hull}_V$
is compact and has a unique topological group structure with
identity $V$ such that $\Z \ni k \mapsto S^k V \in
\mathrm{hull}_V$ is a homomorphism. Moreover, the group structure
is Abelian and there exist arbitrarily small compact open
neighborhoods of $V$ in $\mathrm{hull}_V$ which are finite index
subgroups.
\end{lemma}

In particular, $\Omega = \mathrm{hull}_V$ is a Cantor group, $T =
S|_\Omega$ is a minimal translation, and every element of $\Omega$
may be written in the form \eqref{pot} with the continuous
function $f(\omega) = \omega(0)$.

The second lemma (see \cite[Lemma~2.2]{a}) addresses the converse:

\begin{lemma}
Suppose $\Omega$ is a Cantor group, $T : \Omega \to \Omega$ is a
minimal translation, and $f \in C(\Omega,\R)$. Then, for every
$\omega \in \Omega$, the element $V_\omega$ of $\ell^\infty(\Z)$
defined by $V_\omega (n) = f(T^n \omega)$ is limit-periodic and we
have $\mathrm{hull}_{V_\omega} = \{ V_{\tilde \omega} \}_{\tilde
\omega \in \Omega}$.
\end{lemma}

These two lemmas show that a study of limit-periodic potentials
can be carried out by considering potentials of the form
\eqref{pot} with a minimal translation $T$ of a Cantor group
$\Omega$ and a continuous sampling function $f$. As shown for the
first time in the context of limit-periodic potentials by Avila in
\cite{a}, it is often advantageous to fix $\Omega$ and $T$ and to
vary $f$. This is what we will do in this paper and also in the
forthcoming paper \cite{dg}.

So fix a minimal translation $T$ of a Cantor group $\Omega$. Next,
we discuss the set of periodic sequences that arise in the
representation \eqref{pot}.

\begin{definition}
A sampling function $f \in C(\Omega,\R)$ is called periodic if
there is a period $p \in \N$ such that for every $\omega \in
\Omega$, we have $f(T^p \omega) = f(\omega)$.
\end{definition}

The next lemma (see \cite[Subsubsection~2.3.2]{a}) shows how to
construct periodic sampling functions:

\begin{lemma}
For every $f \in C(\Omega,\R)$ and every compact subgroup
$\Omega_0 \subset \Omega$ of finite index,
$$
f_{\Omega_0}(\omega) = \int_{\Omega_0} f(\omega + \tilde \omega)
\, d\mu_{\Omega_0}(\tilde \omega)
$$
is periodic. Here, $\mu_{\Omega_0}$ denotes Haar measure on
$\Omega_0$.
\end{lemma}

Combining this construction with the fact that one can find
compact subgroups of finite index in any neighborhood of the
identity, we arrive at the following way of looking at the
embedding of periodic sampling functions into $C(\Omega,\R)$;
compare \cite[Section~3]{a}:

\begin{lemma}
There exists a decreasing sequence of Cantor subgroups $\Omega_k
\subset \Omega$ with finite index $p_k$ such that $\bigcap
\Omega_k = \{ 0 \}$. Let $P_k$ be the sampling functions which are
defined on $\Omega / \Omega_k$. Elements of $P_k$ are periodic
with period $p_k$. Conversely, every periodic sampling function
belongs to some $P_k$. Thus, $P = \bigcup P_k$ is the set of all
periodic sampling functions and it is dense in $C(\Omega,\R)$.
\end{lemma}

We refer the reader to \cite{rz} and \cite{w} for further material
and a comprehensive discussion of Cantor groups.

\section{Some Features of Periodic Potentials}\label{s.periodic}

In this section we summarize some aspects of the spectral theory
of Schr\"odinger operators with periodic potentials. Since
limit-periodic potentials are obtained as uniform limits of
periodic potentials, much of their spectral analysis is based on
approximation by periodic potentials and hence a good
understanding of the periodic theory is essential. We will focus
here on the properties of periodic operators that are of immediate
interest to us and refer the reader to \cite{l,s.sz,teschl,t} for
more details.

Consider a potential $V : \Z \to \R$ that is periodic with period
$p$, that is,
$$
V(n+p) = V(n) \quad \text{for every } n \in \Z.
$$
This gives rise to a periodic Schr\"odinger operator in
$\ell2(\Z)$, given by
\begin{equation}\label{peroper}
[H \psi](n) = \psi(n+1) + \psi(n-1) + V(n) \psi(n).
\end{equation}
We will link the spectral properties of $H$ to properties of the
solutions of the associated difference equation
\begin{equation}\label{pereve}
u(n+1) + u(n-1) + V(n)u(n) = Eu(n).
\end{equation}
These solutions are generated by the so-called transfer matrices
$$
T_n^{(E,V)} = S^{(E,V)}_{n-1} \dots S^{(E,V)}_0 ,
$$
where
$$
S^{(E,V)}_i = \begin{pmatrix} E-V(i)& -1\\
1& 0 \end{pmatrix}.
$$
The transfer matrix over one period, $T_p^{(E,V)}$, plays a
special role as we will see below.

For $k\in \mathbb{R}$ and $l \in \mathbb{Z}$, we define
$$
A_l^{k} = \begin{pmatrix} V(l)&1&\ &\ &\ &e^{-ikp}\\1&V(l+1)&1\\\
&1&V(l+2)&1\\\ &\ &\ddots&\ddots&\ddots \\\ &\ &\ &\ &\ &1\\e^{ikp}&\ &\ &\ &1&V(l+p-1)
\end{pmatrix} .
$$
These matrices are the restrictions of $H$ to intervals of length
$p$ with suitable self-adjoint boundary conditions. The importance
of this choice of boundary condition lies in its connection to the
existence of special solutions of \eqref{pereve}.

The following proposition summarizes several important results
concerning the operator, the difference equation, and the
matrices; see the references mentioned above for proofs. The
results listed here are usually referred to as Floquet-Bloch
theory.

\begin{prop}\label{p.periodicbasics}
{\rm (a)} We have $E \in \sigma(H)$ if and only if
\eqref{pereve} has a solution $\{ u(n)\}$
obeying
$$
u(n+p)=e^{ikp}u(n)
$$
for all $n$ and some real number $k$. In this case, $\tilde{u} =
\left<u(n)\right>_{n=l}^{l+p-1}$ is an eigenvector of the matrix
$A_l^k$ corresponding to eigenvalue $E$.

{\rm (b)} The $p$ eigenvalues of $A_l^{k}$ are independent of $l$
and
$$
\sigma(H) = \bigcup_k \sigma(A_l^{k}).
$$

{\rm (c)} The characteristic polynomial of $A_l(k)$ obeys
$$
\det (E - A_l(k)) = \Delta(E)- 2\cos kp,
$$
where $\Delta (E) = \mathrm{Tr} \, T_p^{(E,V)}$. We have
$$
\sigma(H) = \{ E : |\Delta (E)|\le 2 \}.
$$
The set $\sigma(H)$ is made of $p$ bands such that on each band,
$\Delta(E)$ is either strictly increasing or strictly decreasing.

{\rm (d)} If $E$ is in the boundary of some band, we have
$\Delta(E) = \pm 2$. Moreover, if two different bands intersect,
then their common boundary point satisfies $T_p^{(E,V)} = \pm I$.
\end{prop}

The function $\Delta$ is called the discriminant associated with
the periodic potential~$V$. It is a polynomial of degree $p$ with
real coefficients.

The other important consequence of periodicity is the existence of
a direct integral decomposition. This will be described next.

As we have seen above, we can treat  $E \in \sigma(H)$ as a
function of the variable $k \in [0,\frac{\pi}{p}]$. For each band,
the association $k \mapsto E$ is one-to-one and onto. Moreover, if
we consider energies in the interior of a band, that is, with
$\Delta(E) \in (-2,2)$ or $k \in (0,\frac{\pi}{p})$, then there
are linearly independent solutions $\varphi^\pm(E)$ of
\eqref{pereve} with
$$
\varphi^\pm_{n+lp}(E) = e^{\pm ilkp} \varphi^\pm_n (E).
$$
It is easy to see that one can normalize these solutions by
requiring
$$
\varphi^\pm_0(E) > 0
$$
and
$$
\sum_{j = 0}^{p-1} \left|\varphi^\pm_{j}(E)\right|^2 = 1.
$$
With this normalization, we have
$$
\varphi^-(E) = \overline{\varphi^+(E)}.
$$

Next, we define for $u = \{u_n\}_{n \in \Z}$ of finite support,
$$
\hat u^\pm (E) = \sum_{n \in \Z} \overline{\varphi^\pm_n(E)} u_n.
$$
We also define the measure $d\rho$ on $\sigma(H)$ by
$$
d\rho(E) = \frac{1}{\pi} \left|\frac{dk}{dE} (E) \right| \, dE.
$$

Then, we have the following result; see the references listed
above for a proof.

\begin{prop}\label{p.dirint}
The map $u \mapsto \hat u$ extends to a unitary map from
$\ell2(\mathbb{Z})$ to $L2(\sigma(H), d\rho ; \mathbb{C}2)$.
Its inverse is given by
$$
\left( {\check{f}} \right)_n = \frac{1}{2} \int_{\sigma(H)} \left[
\varphi^+_n(E) f^+(E) + \varphi^-_n(E) f^-(E) \right] \, d\rho(E)
$$
Moreover, we have that
$$
\widehat{Hu}^\pm (E) = E \hat{u}^\pm (E).
$$
\end{prop}

Here, we use $f^\pm(E)$ for the two components of a $\C2$-valued
function $f \in L2(\sigma(H), d\rho ; \mathbb{C}2)$.

Proposition~\ref{p.dirint} shows that $H$ has purely absolutely
continuous spectrum (of multiplicity two). More precisely, the
spectral measure associated with the operator $H$ with periodic
potential $V$ and a finitely supported $u \in \ell2(\Z)$ is given
by
\begin{equation}\label{mu}
d\mu_{V,u}(E) = g_{V,u}(E) \, dE
\end{equation}
with density
\begin{equation}\label{G}
g_{V,u}(E)= \frac{1}{2\pi} \left( |\hat{u}^+(E)|^2 +
|\hat{u}^-(E)|^2 \right) \left|\frac{dk}{dE} (E) \right|
\end{equation}
for $E \in \sigma(H)$ (we set $g_{V,u}(E)$ equal to zero outside
of $\sigma(H)$).

Let us derive some consequences of these results.

\begin{lemma}\label{l.t12}
For every $t \in (1,2)$, there exists a constant $D =
D(\|V\|_\infty,p,t)$ such that
\begin{equation}\label{rho}
\int_{\sigma(H)} \left| \frac{dk}{dE}(E) \right|^t \, dE \le D.
\end{equation}
\end{lemma}

\begin{proof}
By Proposition~\ref{p.periodicbasics}, we have
$$
\left| \frac{dk}{dE}(E) \right| = \left| \frac{\Delta' (E)}{2p \,
\sin(kp)} \right|.
$$
Since we can bound $|\Delta' (E)|$ by a
$(\|V\|_\infty,p)$-dependent constant and $\int_0^{\pi}
(\sin(x))^{1-t} \, dx < \infty$, we have the following estimates,
$$
\int_{\sigma(H)} \left| \frac{dk}{dE}(E) \right|^t \, dE \lesssim
\int_0^{\frac{\pi}{p}} \left| \frac{1}{2p \, \sin(kp)}
\right|^{t-1} \, dk \lesssim \int_0^{\frac{\pi}{p}}
|\sin(kp)|^{1-t} \, dk
$$
and the last integral may be bounded by a $t$-dependent constant.
\end{proof}

\begin{lemma}\label{l.one}
Let $u \in \ell2(\mathbb{Z})$ have finite support. Then, for
every $t \in (1,2)$, there exists a constant $Q =
Q(\|V\|_\infty,p,u,t)$ such that
\begin{equation}\label{Gexpre}
\int_{\sigma(H)} \left| g_{V,u}(E) \right|^t \, dE \le Q.
\end{equation}
\end{lemma}

\begin{proof}
Since $u$ has a finite support, we can find a constant $M =
M(p,u)$ such that $|\hat{u}^\pm(E)|^2 \leq M$. Thus, by \eqref{G}
we have
\begin{align*}
\int_{\sigma(H)} \left| g_{V,u} (E) \right|^t \, dE & =
\int_{\sigma(H)} \left[ \frac{1}{2\pi} \left( |\hat{u}^+(E)|^2 +
|\hat{u}^-(E)|^2 \right) \left| \frac{dk}{dE}(E) \right| \right]^t \, dE \\
& \le \left[ \frac{M}{\pi} \right]^t \int_{\sigma(H)} \left|
\frac{dk}{dE}(E) \right|^t
\, dE \\
& \le \left[ \frac{M}{\pi} \right]^t D
\end{align*}
with the constant $D$ from Lemma~\ref{l.t12}.
\end{proof}

\begin{lemma}\label{l.two}
Let $(X,d\mu)$ be a finite measure space, let $r>1$ and let $f_n,f
\in L^r$ with $\sup_n\| f_n \|_r < \infty$. Suppose that $f_n(x)
\rightarrow f(x)$ pointwise almost everywhere. Then, $\| f_n-f
\|_p \rightarrow 0$ for every $p < r$.
\end{lemma}

\begin{proof}
This is \cite[Lemma~2.6]{as}.
\end{proof}

\begin{lemma}\label{l.three}
Suppose $u \in \ell2(\mathbb{Z})$ has finite support and $V_n,V :
\Z \to \R$ are $p$-periodic and such that $\| V_n - V \|_\infty
\rightarrow 0$ as $n \rightarrow \infty$. Then, for any $t \in
(1,2)$, we have
$$
\int_\R \left|g_{V_n,u}(E) - g_{V,u}(E) \right|^t \, dE
\rightarrow 0
$$
as $n \to \infty$.
\end{lemma}

\begin{proof}
By Lemmas~\ref{l.one} and \ref{l.two} we only need to prove
pointwise convergence. Given the explicit identity \eqref{G},
pointwise convergence follows readily from the following two
facts: the discriminant of the approximants converges pointwise to
the discriminant of the limit and the matrices $A_l^k$ associated
with the approximants converge pointwise to those associated with
the limit and therefore so do the associated eigenvectors.
\end{proof}

\section{Cantor Spectrum}

\begin{theorem}\label{t.cantor}
Let $\Omega$ be a Cantor group and let $T: \Omega \rightarrow
\Omega$ be a minimal translation. Then there exists a dense
$G_{\delta}$ set $\mathcal{C} \subseteq C(\Omega,\mathbb{R})$ such
that for every $f \in \mathcal{C}$ and $\omega \in \Omega$, the
spectrum of the operator $H_\omega$ given by \eqref{oper} is a
Cantor set.
\end{theorem}

Fix $\Omega$ and $T$ as in the theorem throughout this section. By
minimality, for given $f \in C(\Omega,\R)$, the spectrum of
$H_\omega$ is independent of $\omega$. For notational convenience,
we will denote this set by $\Sigma(f)$. Let $P = \bigcup P_k$ be
the set of periodic sampling functions in $C(\Omega,\R)$ and
denote the associated periods by $p_k$.

\begin{lemma}\label{l.gdelta}
$\mathfrak{N} := \{f \in C(\Omega,\mathbb{R}) : \Sigma(f)\ has\
empty\ interior \}$ is a $G_\delta$ set.
\end{lemma}
\begin{proof}
This is essentially \cite[Lemma~1.1]{as}.
\end{proof}

\begin{lemma}\label{l.opengaps}
For every $f \in P_k$, $k \in \N$ and every $\varepsilon > 0$,
there exists $\tilde{f}$ in $P_k$ satisfying $\| f - \tilde{f} \|
< \varepsilon$ such that $\Sigma(\tilde{f})$ has exactly $p_k$
components, that is, its $p_k - 1$ gaps are all open.
\end{lemma}

\begin{proof}
This follows from the proof of \cite[Claim 3.4]{a}. For the
reader's convenience, we provide a proof. Let $f \in P_k$ be
given. By $\omega$-independence of the spectrum, we may choose and
fix an arbitrary $\omega \in \Omega$ for the purpose of this
proof. Next, given $\varepsilon > 0$, let $M$ be large enough so
that $\frac{2p_k + 1}{M} < \varepsilon$. Then, for $1 \leq t \leq
2p_k+1$, there is $\tilde{f}_t \in P_k$ with
$$
\tilde{f}_t(T^i \omega) = f(T^i \omega),\; 0\leq\ i\ \leq p_k-2\;
\text{ and } \; \tilde{f}_t(T^{p_k-1} \omega) = f(T^{p_k-1}
\omega) + \frac{t}{M}.
$$
Obviously, $\| \tilde{f}_t - f \| < \varepsilon$ and we claim that
there exists some $t$ in this range such that the spectrum
associated with $\tilde{f}_{t}$ has exactly $p_k$ components.

Suppose this claim fails. Then, for every $t$ in this range, there
exists by Proposition~\ref{p.periodicbasics}.(d) an energy $E_t
\in \Sigma (\tilde{f}_t)$ with $T_{p_k}^{(E_t,\tilde{f}_t)} = \pm
id$. That is,
$$
\begin{pmatrix} E_t-f(T^{p_k-1} \omega)-\frac{t}{M}&-1\\1&0 \end{pmatrix}
\begin{pmatrix} E_t-f(T^{p_k-2} \omega)&-1\\1&0 \end{pmatrix}
\cdots
\begin{pmatrix} E_t-f(\omega)&-1\\1&0 \end{pmatrix}
= \pm id.
$$
Since $f$ is $p_k$-periodic, we get from this
\begin{equation}\label{e.relation}
T_{p_k}^{(E_t,f)} = \pm \begin{pmatrix}1&\frac{t}{M}\\0&1
\end{pmatrix}.
\end{equation}
Indeed, rewriting the identity above, we find that
\begin{align*}
T_{p_k}^{(E_t,f)} & = \pm id + \begin{pmatrix} \frac{t}{M} & 0
\\ 0 & 0 \end{pmatrix} \begin{pmatrix} E_t-f(T^{p_k-2} \omega)&-1\\1&0 \end{pmatrix}
\cdots
\begin{pmatrix} E_t-f(\omega)&-1\\1&0 \end{pmatrix} \\
& = \pm id + \begin{pmatrix} \frac{t}{M} & 0
\\ 0 & 0 \end{pmatrix} \begin{pmatrix} E_t-f(T^{p_k-1} \omega)&-1\\1&0
\end{pmatrix}^{-1} T_{p_k}^{(E_t,f)}
\end{align*}
so that
\begin{align*}
T_{p_k}^{(E_t,f)} & = \pm \left( id - \begin{pmatrix} \frac{t}{M}
& 0 \\ 0 & 0 \end{pmatrix} \begin{pmatrix} E_t-f(T^{p_k-1}
\omega)&-1\\1&0
\end{pmatrix}^{-1} \right)^{-1} \\
& = \pm \left( id - \begin{pmatrix} \frac{t}{M} & 0 \\ 0 & 0
\end{pmatrix} \begin{pmatrix} 0 & 1 \\ -1 & E_t-f(T^{p_k-1} \omega)
\end{pmatrix} \right)^{-1} \\
& = \pm \left( id - \begin{pmatrix} 0 & \frac{t}{M} \\ 0 & 0
\end{pmatrix} \right)^{-1} \\
& = \pm \begin{pmatrix} 1 & - \frac{t}{M} \\ 0 & 1
\end{pmatrix}^{-1} \\
& = \pm \begin{pmatrix} 1 & \frac{t}{M} \\ 0 & 1
\end{pmatrix}.
\end{align*}

The relation \eqref{e.relation} implies that if $t\neq t'$, we
have $T_{p_k}^{(E_t,f)} \neq T_{p_k}^{(E_{t'},f)}$, and therefore
$E_t \neq E_{t'}$. But there are at most $2 p_k$ values of $E$ for
which $tr\ T_{p_k}^{(E,f)} = \pm 2$; contradiction.
\end{proof}

\begin{lemma}\label{l.measure}
Let $f \in P$ have period $p$. Then we have
\\
{\rm (1)} The measure of each band of $\Sigma(f)$ is at most
$\frac{2\pi}{p}$.
\\
{\rm (2)} Let $C\ge 1$ be such that for every $E \in \Sigma(f)$,
there exist $\omega \in \Omega$ and $k\ge 1$ such that $\|
T_k^{(E,V_\omega)} \| \ge C$. Then, the total measure of
$\Sigma(f)$ is at most $\frac{4\pi p}{C}$.
\end{lemma}

\begin{proof}
This is \cite[Lemma~2.4]{a}.
\end{proof}

\begin{proof}[Proof of Theorem~\ref{t.cantor}.] The proof is close in spirit to the
proof of \cite[Theorem~1]{as}. By Lemma~\ref{l.gdelta}, we only
need to prove that $\mathfrak{N}$ is dense. Since $P = \bigcup
P_k$ is dense in $C(\Omega,\R)$, it suffices to show that, given
$f\in P$ and $\varepsilon > 0$, there is a potential $\tilde{f}$
such that $\| \tilde{f}-f \| < \varepsilon$ and
$\Sigma(\tilde{f})$ is nowhere dense.

So let $f \in P$ and $\varepsilon > 0$ be given. Write $f$ as $f =
\sum^N_{j=0} a_j W_j$ with $W_j \in P_j$. We construct $s_0 =
\sum^{N}_{i=0} a_n^{(0)} W_i$ so that $\| s_0\| < \varepsilon$ and
$f_0 = f+s_0$ has all $n_N-1$ gaps open. (This is possible due to
Lemma~\ref{l.opengaps}.)

Suppose $s_0,s_1, \ldots ,s_{k-1}$ are picked. Let $\alpha_{k-1}$
be the minimal gap size of $f_{k-1}$ and define $\beta_k = \min \{
\alpha_0, \alpha_1, \ldots ,\alpha_{k-1} \}$. Applying
Lemma~\ref{l.opengaps}, we pick $s_k = \sum^{N+k}_{i=0} a_i^{(k)}
W_i$ so that
\begin{align}
\label{e.choice4} \| s_k \| & < \frac{\varepsilon}{2^k} , \\
\label{e.choice5} \| s_k \| & < \frac{1}{3}\frac{\beta_k}{2^k} , \\
\label{e.choice6} f_k & = f + \sum^k_{j=0}s_j \text{ has all gaps
open.}
\end{align}

The limit of $f_k$ exists by \eqref{e.choice4}, let $\tilde{f} =
\lim_{k \rightarrow \infty}f_k$. By construction, we have $\|
\tilde f - f \| < \varepsilon$. We claim that $\Sigma(\tilde f)$
is nowhere dense; equivalently, its complement is dense.

Given $E \in \Sigma(\tilde{f})$ and $\tilde \varepsilon
> > 0$, we can pick $k$ large enough so that
\begin{align}
\label{e.choice1} \| \tilde{f}-f_k\| & < \frac{\tilde
\varepsilon}{3} , \\
\label{e.choice2} \frac{2\pi}{n_{N+k}} & < \frac{\tilde
\varepsilon}{3} , \\
\label{e.choice3} \frac{\varepsilon}{2^{k-1}} & < \frac{\tilde
\varepsilon}{3} .
\end{align}
By \eqref{e.choice1}, there exists $E' \in \Sigma(f_k)$ such that
$|E' - E| < \frac{\tilde \varepsilon}{3}$. Moreover, by
Lemma~\ref{l.measure} and \eqref{e.choice2}, we can find
$\tilde{E}$ in a gap of $\Sigma(f_k)$ such that $|E' - \tilde{E}|
< \frac{\tilde \varepsilon}{3}$. Write this gap of $\Sigma(f_k)$
that contains $\tilde E$ as $(a-\delta,a+\delta)$. By definition,
we have $2\delta \ge \beta_{k+1}$. By \eqref{e.choice5},
$$
\| \tilde{f}-f_k \| = \left\| \sum^{\infty}_{j = k+1} s_j \right\|
< \frac{\beta_{k+1}}{3} \left( \frac{1}{2^{k+1}} +
\frac{1}{2^{k+2}} + \cdots \right) \leq \frac{\delta}{3},
$$
so we have $(a-\frac{\delta}{3},a+\frac{\delta}{3})\cap
\Sigma(\tilde f)= \emptyset$.

We claim that there exists $\delta' \in [\frac{\delta}{3} ,
\delta)$ such that $(a - \delta' , a + \delta') \cap
\Sigma(\tilde{f}) = \emptyset$ and $|\delta' - \delta| <
\frac{\varepsilon}{2^k}$. As we saw above, we may arrange for
$\delta' \ge \frac{\delta}{3}$. Suppose that it is impossible to
find such a $\delta'$ with $|\delta' - \delta| <
\frac{\varepsilon}{2^k}$. Then, there will be a point $x \in
\Sigma(\tilde{f})$ such that $[x - \frac{\varepsilon}{2^k} , x +
\frac{\varepsilon}{2^k}] \subseteq (a-\delta,a+\delta)$. Then we
get a contradiction to the already established fact
$[x-\frac{\varepsilon}{2^k},x+\frac{\varepsilon}{2^k}] \cap
\Sigma(f_k) = \emptyset$ since \eqref{e.choice4} implies
$$
\| \tilde{f} - f_k \| = \left\| \sum^{\infty}_{j = k+1} s_j
\right\| < \frac{\varepsilon}{2^k}.
$$

We can choose an energy $\hat{E}$ in the gap of
$\Sigma(\tilde{f})$ that contains $(a-\delta',a+\delta')$ such
that $|\hat{E} - \tilde{E}|\leq \frac{\varepsilon}{2^k} <
\frac{\tilde \varepsilon}{3}$, where the second inequality follows
from \eqref{e.choice3}. Moreover, since we also have
$|\tilde{E}-E'| < \frac{\tilde \varepsilon}{3}$ and $|E' - E| <
\frac{\tilde \varepsilon}{3}$, it follows that $|\hat{E}-E| <
\tilde \varepsilon $. This shows that $\R \setminus
\Sigma(\tilde{f})$ is dense and completes the proof.
\end{proof}

\noindent\textit{Remark.} The reader may notice that the statement of Theorem~\ref{t.cantor} is a consequence of \cite[Corollary~1.2]{a}. However, the main purpose of this section is the method of proof presented here, which is direct and flexible enough so that it can be used to also ensure absolutely continuous spectrum (this will be done in the next section). In particular, the Cantor spectra constructed here may have positive Lebesgue measure, whereas the Cantor spectra generated in the proof of \cite[Corollary~1.2]{a} always have zero Lebesgue measure.

\section{Absolutely Continuous Spectrum}

\begin{theorem}\label{t.ac}
Let $T: \Omega \rightarrow  \Omega$ be a minimal translation of a
Cantor group. Then there exists a dense set of $f \in
C(\Omega,\mathbb{R})$ so that for every $\omega \in \Omega$, the
spectrum of $H_\omega$ is a Cantor set of positive Lebesgue
measure and all spectral measures are purely absolutely
continuous.
\end{theorem}

\begin{proof}
The idea is to modify the construction from the proof of
Theorem~\ref{t.cantor}. Thus, we will again start with an
arbitrarily small ball in $C(\Omega,\R)$ and construct a point in
this ball for which the associated Schr\"odinger operator has both
Cantor spectrum and purely absolutely continuous spectrum. The
presence of absolutely continuous spectrum then also implies that
the Lebesgue measure of the spectrum is positive.

Fix $t\in (1,2)$ and let $u \in \ell2(\mathbb{Z})$ have finite
support. In going through the construction in the proof of
Theorem~\ref{t.cantor}, pick $s_k$ so that in addition to the
conditions above, we have
\begin{equation}\label{Gfinal}
\left(\int_{-\infty}^{\infty} \left|g_u^{k-1}(E) -
g_u^{k}(E)\right|^t\right)^{\frac{1}{t}} \leq \frac{1}{2^k},
\end{equation}
where $g_u^{k}$ is the density of the spectral measure associated
with $u$ and the periodic potential $n \mapsto f_k(T^n \omega)$,
with the estimate above being uniform in $\omega \in \Omega$. This
is possible due to Lemma~\ref{l.three}.

By Lemma~\ref{l.one}, there exists a constant $Q(u,t) < \infty$
such that $\int_\R \left|g_u^k(E)\right|^t \, dE \leq Q(u,t)$.

Now fix any $\omega \in \Omega$. Let $A$ be a finite union of open
sets. If $P_A^k$ is the spectral projection for the potential $n
\mapsto f_k(T^n \omega)$ and $P_A$ is the spectral projection for
the potential $n \mapsto \tilde f(T^n \omega)$, it follows that
$\langle u,P_A u \rangle \le \limsup_{k\rightarrow \infty} \langle
u, P_A^{f_k} u \rangle$ since $\|f_k - \tilde{f}\|_\infty \to 0$
and hence the associated Schr\"odinger operators converge in norm.

Applying H\"older's inequality, we find
$$
\langle u, P_A u \rangle \le \limsup_{k \to \infty} \int_A
g_u^k(E) \, dE \le Q(u,t) |A|^{\frac{1}{q}},
$$
where $\frac{1}{q}+\frac{1}{t}=1$ and $| \cdot |$ denotes Lebesgue
measure. This shows that the spectral measure associated with $u$
and the Schr\"odinger operator with potential $n \mapsto \tilde
f(T^n \omega)$ is absolutely continuous with respect to Lebesgue
measure. Since this holds for every finitely supported $u$, it
follows that this operator has purely absolutely continuous
spectrum.
\end{proof}

Notice that this establishes part (a) of Theorem~\ref{t.main}.

\section{Absence of Point Spectrum}

\begin{definition}
A bounded function $V : \Z \rightarrow \R$ is called a Gordon
potential if there are positive integers $q_k \rightarrow \infty$
such that
$$
\max_{1\leq n \leq q_k} \left| V(n) - V(n \pm q_k) \right| \leq
k^{-q_k}
$$
for every $k \ge 1$.
\end{definition}

Clearly, if $V$ is a Gordon potential, then so is $\lambda V$ for
every $\lambda \in \R.$

\begin{lemma}
Suppose $V$ is a Gordon potential. Then, the operator $H$ given by
\eqref{peroper} has empty point spectrum.
\end{lemma}

This is essentially due to Gordon \cite{g}; see \cite{dp} for the
modification of the argument necessary to prove the result as
stated.

\begin{theorem}\label{t.absencepp}
Let $\Omega$ be a Cantor group and let $T: \Omega \rightarrow
\Omega$ be a minimal translation. Then there exists a dense
$G_\delta$-set $\mathcal{G} \subseteq C(\Omega,\mathbb{R})$ such
that for every $f \in \mathcal{G}$ and $\omega \in \Omega$, the
potential $V_\omega$ given by \eqref{pot} is a Gordon potential.
\end{theorem}

\begin{proof}
We know the set of periodic potentials is dense in
$C(\Omega,\mathbb{R})$. For $j,k \in \mathbb{N}$, let
$$
\mathcal{G}_{j,k} = \left\{ f \in C(\Omega,\mathbb{R}) : \text{
there is a $j$-periodic $f_j$ such that } \|f - f_j\| <
{\frac{1}{2}}(jk)^{(-jk)} \right\}.
$$
Clearly, $\mathcal{G}_{j,k}$ is open. For $k \in \N$, let
$$
\mathcal{G}_k = \bigcup_{j = 1}^\infty \mathcal{G}_{j,k}.
$$
The set $\mathcal{G}_k$ is open by construction and dense since it
contains all periodic sampling functions. Thus,
$$
\mathcal{G} = \bigcap_{k=1}^\infty \mathcal{G}_k
$$
is a dense $G_\delta$ subset of $C(\Omega,\mathbb{R})$. We claim
that for every $f \in \mathcal{G}$ and every $\omega \in \Omega$,
the potential $V_\omega$ given by \eqref{pot} is a Gordon
potential.

So let $f \in \mathcal{G}$ and $\omega \in \Omega$ be given. Since
$f \in A_k$ for every $k \in \mathbb{N}$, we can find
$j_k$-periodic $f_{j_k}$ satisfying
$$
\| f - f_{j_k} \| < {\frac{1}{2}} (j_k k)^{- j_k k}.
$$
Let $q_k = j_k k$, so that $q_k \rightarrow \infty$ as $k
\rightarrow \infty$. Then, we have
\begin{align*}
\max_{1\leq n \leq q_k} & \|V_\omega (n) - V_\omega (n \pm q_k)\|
= \max_{1\leq n \leq q_k} \|f(T^n \omega) - f(T^{n \pm
q_k} \omega)\| \\
& = \max_{1 \leq n
\leq j_k k} \|f(T^n \omega) - f_{j_k}(T^n \omega) + f_{j_k}(T^{n \pm j_k k} \omega) - f(T^{n \pm j_k k} \omega)\| \\
& \leq \max_{1 \leq n \leq j_k k} \|f(T^n \omega) - f_{j_k}(T^n
\omega)\| + \max_{1 \leq n \leq j_k k} \|f(T^{n \pm j_k k} \omega)
- f_{j_k}(T^{n \pm j_k k} \omega)\| \\
& < {\frac{1}{2}}(j_k k)^{-j_k k} + {\frac{1}{2}}(j_k k)^{-j_k k} \\
& \le k^{-j_k k} \\
& = k^{-q_k}.
\end{align*}
It follows that $V_\omega$ is a Gordon potential.
\end{proof}

\begin{lemma}\label{l.avilazeroleb}
Let $T: \Omega\rightarrow \Omega $ be a minimal translation of a
Cantor group. For a dense $G_\delta$ set of $f \in C(\Omega,\R)$,
and for every $\lambda \neq 0$, the Schr\"{o}dinger operator with
potential $\lambda f(T^n \omega)$ has a spectrum of zero Lebesgue
measure for every $\omega \in \Omega$.
\end{lemma}

\begin{proof}
This is \cite[Corollary 1.2.]{a}.
\end{proof}

We can now finish the proof of our main theorem.

\begin{proof}[Proof of Theorem~\ref{t.main}.{\rm (}b{\rm )}.]
Since the intersection of two dense $G_\delta$ sets is again a
dense $G_\delta$ set and zero-measure spectrum precludes
absolutely continuous spectrum, the result follows directly from
Theorem~\ref{t.absencepp} and Lemma~\ref{l.avilazeroleb}.
\end{proof}

\end{document}